\documentclass[11pt]{article}

\usepackage{tikz}
\usepackage{caption}
\usepackage{subcaption}

\usepackage{times}
\usepackage{amsmath}
\usepackage{amsfonts}
\usepackage{amssymb}

\usepackage{tabularx}
\usepackage[linkcolor=black,urlcolor=black,citecolor=black]{hyperref}
\usepackage[all]{hypcap}
\hypersetup{linktocpage,colorlinks=true,pdfborder={0 0 0}}

\usepackage{geometry}
\title{Approximate Solutions of 4-regular Matchstick Graphs\\with 50 -- 62 Vertices\\[5mm]\Large{A catalog of examples which contain 2 or 3 forbidden distances}}

\author{Mike Winkler}

\date{Fakult\"at f\"ur Mathematik\\ Ruhr-Universit\"at Bochum, Germany\\ mike.winkler@ruhr-uni-bochum.de\\www.mikematics.de\\[5mm]June 14, 2020\\[5mm]}

\begin{document}
  
  \maketitle
  
  \begin{abstract}
    A 4-regular matchstick graph is a planar unit-distance graph whose vertices have all degree 4. Examples of 4-regular matchstick graphs are currently known for all number of vertices $\geq$ 52 except for 53, 55, 56, 58, 59, 61, and 62. In this article we present 38 different examples with 50 -- 62 vertices which contain two, three, or four distances which differ slightly from the unit length. These graphs should show why this subject is so extraordinarily difficult to deal with and should also be an incentive for the interested reader to find solutions for the missing numbers of vertices.
  \end{abstract}
  
  \begin{center}\end{center}
  
  \begin{center} 
  	\begin{minipage}{\linewidth}
  		\centering
  		\begin{tikzpicture}[draw=black,font=\sffamily\tiny,
  		myscale/.style={scale=#1, node distance=#1*5em},myscale=1]
  		\begin{scope}[xshift=4cm, myscale=1.2]
  		
  		\foreach \i/\x/\y in {
  			1/1.92462189062386856975/5.78020978657798156775,
  			2/1.10441825531107173930/5.20813811570126983241,
  			3/2.00994867273211230696/4.78385676668239501197,
  			4/1.18974503741931525447/4.21178509580568416482,
  			5/0.28421461999827474232/4.63606644482455809708,
  			6/2.25013346613815423680/4.83467173117978266106,
  			7/2.90623765460078953282/5.58934205250015025968,
  			8/3.23174923011507519988/4.64380399710195135299,
  			9/3.88785341857771049590/5.39847431842231806343,
  			10/4.21336499409199660704/4.45293626302411915674,
  			11/4.86946918255463145897/5.20760658434448675536,
  			12/4.41094248177624681517/4.31892595291759917586,
  			13/5.40982583483231760368/4.36617049744349561280,
  			14/4.95129913405393384807/3.47748986601660980966,
  			15/5.95018248711000463658/3.52473441054250624660,
  			16/4.96572445788796024146/3.34911439829946200319,
  			17/5.61004486451439188244/2.58435874215513017305,
  			18/4.62558683529234659915/2.40873872991208681782,
  			19/5.26990724191877912830/1.64398307376775409949,
  			20/4.27329183709048532336/1.72618851676842965936,
  			21/4.70040753753669360293/0.82199153688387704975,
  			22/3.70379213270840024208/0.90419697988455249860,
  			23/4.13090783315460807756/0.00000000000000000000,
  			24/3.63090783315460807756/0.86602540378443870761,
  			25/3.13090783315460807756/0.00000000000000000000,
  			26/2.63090783315460807756/0.86602540378443870761,
  			27/2.13090783315460852165/0.00000000000000049783,
  			28/1.63090783315460852165/0.86602540378443904068,
  			29/1.13090783315460874370/0.00000000000000066377,
  			30/1.56246121280795824404/0.90208740181302449201,
  			31/0.56545391657730426083/0.82477989077534374918,
  			32/0.99700729623065376117/1.72686729258836746403,
  			33/0.00000000000000000000/1.64955978155068683222,
  			34/0.90949931632193636855/2.06526519841451650095,
  			35/0.09473820666609165941/2.64506200264197754990,
  			36/1.00423752298802804184/3.06076741950580721863,
  			37/0.18947641333218323556/3.64056422373326782349,
  			38/1.09897572965411938206/4.05626964059709749222,
  			39/1.99401459246130752234/1.80417480362604831789,
  			40/1.64602607960155755684/2.74167360091469580752,
  			41/2.61848977766767809783/1.02313015505500159996,
  			42/3.25639329570712687811/1.79324644606103600353,
  			43/2.02084608903060392038/3.66877120124144129321,
  			44/2.27050126480792702210/1.96062895234364642505,
  			45/2.31270396552577928162/3.83078844138487495385,
  			46/3.98126642866591495817/3.17349438605641953615,
  			47/3.46606266738423807183/2.31642666156708898484,
  			48/3.95241578099786217138/3.43024532149071337273,
  			49/2.64532127423697360769/2.88772655267039501936,
  			50/3.25984773076972356165/4.15159786761722227055,
  			51/2.98142212585365040667/3.19114006920694315284}
  		\coordinate (p-\i) at (\x,\y);
  		
  		\foreach \i/\j in {
  			1/2, 1/3, 
  			2/5, 2/4, 
  			3/2, 3/4, 
  			4/5, 
  			5/37, 5/38, 
  			6/7, 6/1, 
  			7/1, 
  			8/7, 8/6, 
  			9/7, 9/8, 
  			10/9, 10/8, 
  			11/9, 11/10, 
  			12/13, 12/11, 
  			13/11, 
  			14/13, 14/12, 
  			15/13, 15/14, 
  			16/17, 16/15, 
  			17/15, 17/19, 
  			18/19, 18/17, 18/16, 
  			20/21, 20/19, 
  			21/23, 21/19, 
  			22/23, 22/21, 22/20, 
  			24/25, 24/23, 
  			25/23, 
  			26/25, 26/24, 
  			27/25, 27/26, 
  			28/27, 28/26, 
  			29/27, 29/28, 
  			30/31, 30/29, 
  			31/29, 
  			32/31, 32/30, 
  			33/31, 33/32, 
  			34/35, 34/33, 
  			35/33, 
  			36/35, 36/34, 
  			37/35, 37/36, 
  			38/37, 38/36, 
  			39/32, 39/30, 
  			40/34, 40/39, 
  			41/28, 41/39, 41/42, 
  			42/24, 42/22, 
  			43/40, 43/49, 43/38, 43/4, 
  			44/41, 44/42, 44/40, 
  			45/3, 45/50, 45/6, 
  			46/18, 46/16, 
  			47/20, 47/46, 
  			48/14, 48/12, 
  			49/45, 49/47, 49/44, 
  			50/10, 50/48, 
  			51/48, 51/50, 51/47, 51/46}
  		\draw[gray,thin] (p-\i) -- (p-\j);
  		
  		\draw[red,thin] (p-43) -- (p-4);
  		\draw[red,thin] (p-42) -- (p-22);
  		\draw[red,thin] (p-45) -- (p-6);
  		
  		\foreach \i in {1,...,51}
  		\fill[black] (p-\i) circle (1pt);
  		
  		\end{scope}
  		\end{tikzpicture}
  		
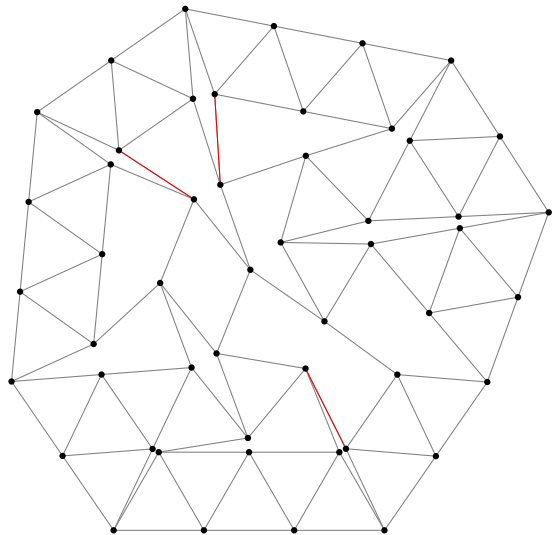
\captionof{figure}{51 vertices}
  	\end{minipage}
  \end{center}
  
  \newpage
  
  \section{Introduction}
  
  A matchstick graph is a planar unit-distance graph. That is a graph drawn with straight edges in the plane such that the edges have unit length, whereby non-adjacent edges do not intersect. We call a matchstick graph 4-regular if every vertex has only degree 4.
  
  Examples of 4-regular matchstick graphs are currently known for all number of vertices $\geq$ 52 except for 53, 55, 56, 58, 59, 61, and 62. For 52, 54, 57, 60, and 64 vertices only one example is known. For a proof we refer the reader to \cite{Existence}. An overview of the currently known examples with 63 -- 70 vertices can be found in \cite{Catalogue}. The currently smallest known example with 52 vertices is the so-called Harborth Graph presented first by Heiko Harborth in 1986 \cite{Harborth}.
  
  \begin{center} 
  	\begin{minipage}{\linewidth}
  		\centering
  		\begin{tikzpicture}[draw=black,font=\sffamily\tiny,
  		myscale/.style={scale=#1, node distance=#1*5em},myscale=1]
  		\begin{scope}[xshift=4cm, myscale=1.0]
  		
  		\foreach \i/\x/\y in {
  			1/0.86128313704930681283/1.17537529692625897226,
  			2/1.78133705523837515550/0.78358353128417235212,
  			3/1.66061171818344721629/1.77626948010835827851,
  			4/2.58066563637251489283/1.38447771446627165837,
  			5/2.70139097342744305408/0.39179176564208617606,
  			6/3.50071955456158256936/0.99268594882418559333,
  			7/3.62144489161651117470/0.00000000000000000000,
  			8/1.42757023375669600540/1.99958335567278444600,
  			9/0.43064156852465335090/2.07789833788345568166,
  			10/0.99692866523204248796/2.90210639662998115540,
  			11/0.00000000000000000000/2.98042137884065239106,
  			12/1.99385733046408497593/2.82379141441930991974,
  			13/2.62333153699726340236/2.04677003549243297087,
  			14/2.98151468716496426126/2.98042137884065638787,
  			15/0.86128313704930337114/4.78546746075504803031,
  			16/1.78133705523836893825/5.17725922639713687090,
  			17/1.66061171818344344153/4.18457327757295161064,
  			18/2.58066563637251134011/4.57636504321503956305,
  			19/2.70139097342743550456/5.56905099203922659967,
  			20/3.50071955456157724029/4.96815680885712751547,
  			21/3.62144489161650362519/5.96084275768131455209,
  			22/1.42757023375669378495/3.96125940200852433293,
  			23/0.43064156852465163006/3.88294441979784998864,
  			24/0.99692866523204237694/3.05873636105132629126,
  			25/1.99385733046408475388/3.13705134326199974737,
  			26/2.62333153699726162600/3.91407272218887891668,
  			27/6.38160664618370265799/4.78546746075505868845,
  			28/5.46155272799463986644/5.17725922639714131179,
  			29/5.58227806504957069222/4.18457327757295427517,
  			30/4.66222414686049457799/4.57636504321504311577,
  			31/4.54149880980557352217/5.56905099203922748785,
  			32/3.74217022867143178644/4.96815680885712751547,
  			33/5.81531954947631657404/3.96125940200853055018,
  			34/6.81224821470835895099/3.88294441979786064678,
  			35/6.24596111800096842614/3.05873636105133650531,
  			36/7.24288978323301613216/2.98042137884066171694,
  			37/5.24903245276892960192/3.13705134326200418826,
  			38/4.61955824623574873300/3.91407272218888246940,
  			39/4.26137509606805053863/2.98042137884065772013,
  			40/6.38160664618371242796/1.17537529692626363520,
  			41/5.46155272799464519551/0.78358353128417601585,
  			42/5.58227806504957069222/1.77626948010836338554,
  			43/4.66222414686050434796/1.38447771446627410086,
  			44/4.54149880980557973942/0.39179176564208650912,
  			45/3.74217022867143711551/0.99268594882418681458,
  			46/5.81531954947631923858/1.99958335567279132938,
  			47/6.81224821470836428006/2.07789833788346367527,
  			48/6.24596111800097286704/2.90210639662998737265,
  			49/5.24903245276893315463/2.82379141441931258427,
  			50/4.61955824623575228571/2.04677003549243519132,
  			51/3.62144489161650495745/1.98537189764837185280,
  			52/3.62144489161650451337/3.97547086003294092293}
  		\coordinate (p-\i) at (\x,\y);
  		
  		\foreach \i/\j in {
  			2/1, 
  			3/1, 3/2, 
  			4/3, 4/2, 
  			5/4, 5/2, 
  			6/4, 6/5, 
  			7/6, 7/5, 7/44, 7/45, 
  			8/1, 
  			9/1, 9/8, 
  			10/9, 10/8, 
  			11/9, 11/10, 11/23, 11/24, 
  			12/10, 12/8, 
  			13/12, 13/3, 13/51, 
  			14/12, 14/13, 14/25, 14/26, 
  			16/15, 
  			17/15, 17/16, 
  			18/16, 18/17, 
  			19/16, 19/18, 
  			20/18, 20/19, 
  			21/19, 21/20, 21/31, 21/32, 
  			22/15, 
  			23/15, 23/22, 
  			24/22, 24/23, 
  			25/22, 25/24, 
  			26/17, 26/25, 26/52, 
  			28/27, 
  			29/27, 29/28, 
  			30/28, 30/29, 
  			31/28, 31/30, 
  			32/30, 32/31, 
  			33/27, 
  			34/27, 34/33, 
  			35/33, 35/34, 
  			36/34, 36/35, 36/47, 36/48, 
  			37/33, 37/35, 
  			38/29, 38/37, 38/52, 
  			39/37, 39/38, 39/49, 39/50, 
  			41/40, 
  			42/40, 42/41, 
  			43/41, 43/42, 
  			44/41, 44/43, 
  			45/43, 45/44, 
  			46/40, 
  			47/40, 47/46, 
  			48/46, 48/47, 
  			49/46, 49/48, 
  			50/42, 50/49, 50/51, 
  			51/6, 51/45, 
  			52/32, 52/20}
  		\draw[gray,thin] (p-\i) -- (p-\j);
  		
  		\foreach \i in {1,...,52}
  		\fill[black] (p-\i) circle (1.3pt);
  		
  		\end{scope}
  		\end{tikzpicture}
  		
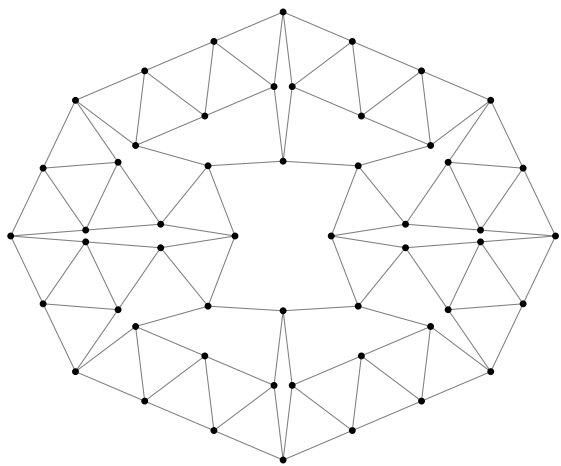
\captionof{figure}{The Harborth Graph}
  	\end{minipage}
  \end{center}
  \quad\\
  In this article we present 38 different examples of 4-regular rigid planar graphs for all number of vertices $\geq$ 50 and $\leq$ 62. These graphs are not matchstick graphs, because each graph contains at least two edges which differ slightly from the unit length. The edges which do not have unit length are colored red.
  
  The geometry of the graphs has been verified with the software \textsc{Matchstick Graphs Calculator} (MGC) \cite{MGC}. This remarkable software created by Stefan Vogel runs directly in web browsers\footnote{For optimal functionality and design please use the Firefox web browser.}. A special version of the MGC contains all graphs from this article and is available on the author's \href{http://mikematics.de/matchstick-graphs-calculator.htm}{website}\footnote{http://mikematics.de/matchstick-graphs-calculator.htm}. The graphs were constructed and first presented by the author between March 21, 2014 and June 5, 2020 in a graph theory internet forum \cite{Thread1} \cite{Thread2}.\footnote{Except Figure 27 by Peter Dinkelacker.}
  
  \newpage
  
  \section{Construction rules and Epsilon graphs}
  
  To get a kind of \textit{fair} approximate solutions we constructed the graphs by using the following rules.
  \begin{itemize}
  	
  	\item The graph must be rigid.
  	
  	\item Equilateral triangles which contain vertices of the outer circle of the graph may not contain forbidden distances. These set of equilateral triangles we denote as the \textit{frame} of the graph.
  	
  	\item The graph may not contain more than three forbidden distances.
  	
  	\item Whenever possible the forbidden distances may not deviate more than 10 percent from the unit length.
  	
  \end{itemize}
  
  The exceptions from the last construction rule apply to the number of vertices for which only one approximate example currently have been found or graphs that we do not want to deprive the interested reader.
  
  Please note that much better approximate solutions are possible if we would ignore the second and/or the third construction rule. The deviation from the unit length of an edge with a forbidden distance becomes smaller the closer these edge is to the outer circle of the graph, or if we distributing the deviation on all edges of the graph.
  
  Without the first three construction rules it is possible to construct flexible graphs whose forbidden distances can be infinitesimally approximated to the unit length. These kind of graphs we will denote as \textit{Epsilon graphs}. Figure 2 shows the smallest possible example with a minimum number of vertices. This Epsilon graph has 27 vertices, a rotational symmetry of order 3, and contains six forbidden distances of equal length.
  \\
  \begin{figure}[!ht]
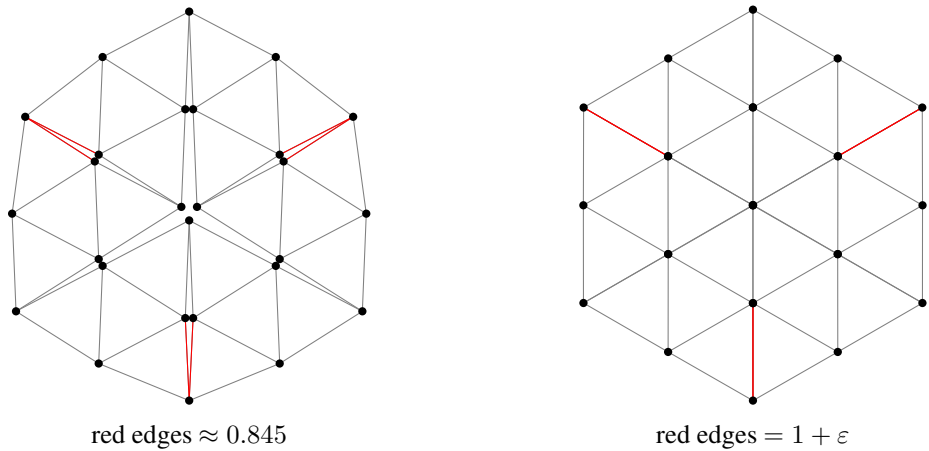
 
  	\centering
  	\begin{minipage}[t]{0.45\linewidth}	
  		\centering		

  		\subcaption*{red edges $\approx0.845$}
  	\end{minipage}
  	\quad
  	\begin{minipage}[t]{0.45\linewidth}	
  		\centering
  		%
	
  		\subcaption*{red edges $=1+\varepsilon$}
  	\end{minipage}
  	\caption{27 vertices}
  \end{figure}
  
  \newpage
  
  Figure 3 and 4 show the smallest known example with only 4 forbidden distances. These Epsilon graph has 42 vertices and a point symmetry.
  \\
  \begin{center}
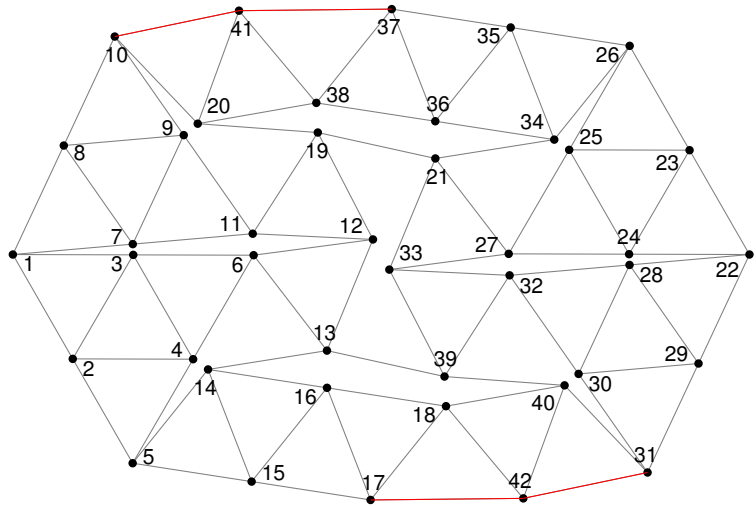

  	\begin{minipage}{\linewidth}
  		\centering	
  		%
  		\captionof{figure}{42 vertices, point symmetry}
  	\end{minipage}
  \end{center}
  {\par\centering$\vert P10,P41\vert=\vert P42,P31\vert\approx1.036$, $\vert P41,P37\vert=\vert P17,P42\vert\approx1.167$\par}
  
  \begin{center}\end{center}
  
  \begin{center}
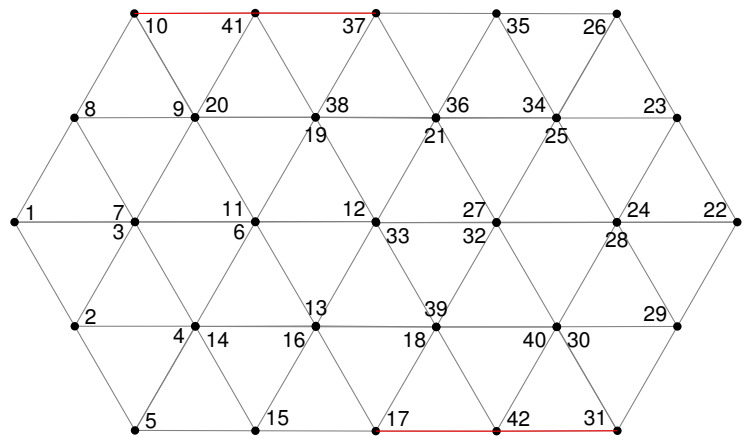

  	\begin{minipage}{\linewidth}
  		\centering	
  		%
  		\captionof{figure}{42 vertices, red edges $=1+\varepsilon$}
  	\end{minipage}
  \end{center}
  
  \quad\\
  
  Even if it contradicts the intuition, Figure 2 (right-sided) and Figure 4 do not show matchstick graphs which contain vertices of degree 3, 4, and 6. They show 4-regular planar graphs which are not matchstick graphs.
  
  \newpage


  \section{Examples of 4-regular planar graphs with 50 -- 62 vertices that look like matchstick graphs}

  \begin{center}
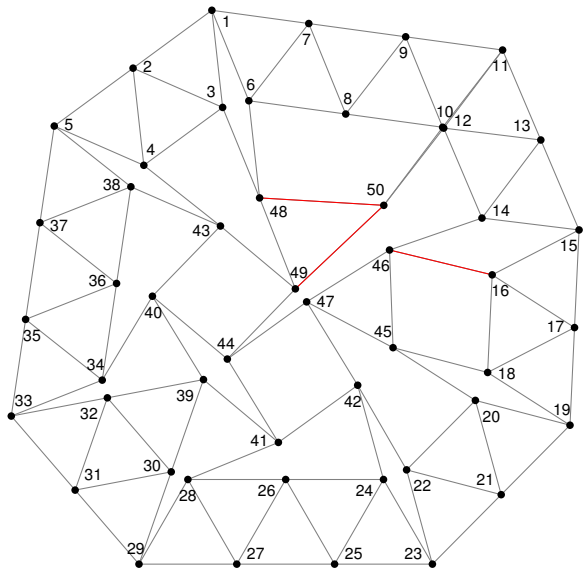
 
  	\begin{minipage}{\linewidth}
  		\centering
  		%
  		\captionof{figure}{50 vertices, asymmetric}
  	\end{minipage}
  \end{center}
  {\par\centering$\vert P46,P16\vert\approx1.0797549592$, $\vert P48,P50\vert\approx1.2721354299$, $\vert P49,P50\vert\approx1.2440648255$\par}
  
  \begin{center}\end{center}

  \begin{center}
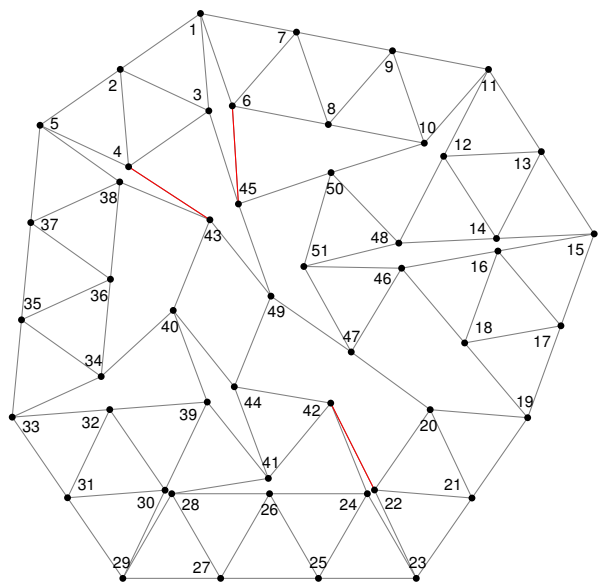
 
  	\begin{minipage}{\linewidth}
  		\centering
  		%
  		\captionof{figure}{51 vertices, asymmetric}
  	\end{minipage}
  \end{center}
  {\par\centering$\vert P42,P22\vert\approx0.99527617909$, $\vert P43,P4\vert\approx0.99277039021$, $\vert P45,P6\vert\approx1.00583136108$\par}

  \begin{center} 
  	\begin{minipage}{\linewidth}
  		\centering
  		
  		%
  		\captionof{figure}{51 vertices, asymmetric}
  	\end{minipage}
  \end{center}
  {\par\centering$\vert P43,38\vert\approx1.0096420153$, $\vert P44,P40\vert\approx0.9924715954$, $\vert P51,P46\vert\approx1.0027811544$\par}

  \begin{center} 
  	\begin{minipage}{\linewidth}
  		\centering
  		
  		%
  		\captionof{figure}{51 vertices, asymmetric}
  	\end{minipage}
  \end{center}
  {\par\centering$\vert P43,P4\vert\approx0.9890758621$, $\vert P45,P6\vert\approx1.0027078452$, $\vert P46,P14\vert\approx0.9922899189$\par}
  
  \begin{center} 
  	\begin{minipage}{\linewidth}
  		\centering
  		
  		%
  		\captionof{figure}{\small$\vert P48,P10\vert\approx1.0000027505$, $\vert P43,P38\vert\approx1.0019150342$, $\vert P42,P22\vert\approx0.9887895316$}
  	\end{minipage}
  \end{center}
  {\par\centering  \par}

  \begin{center} 
  	\begin{minipage}{\linewidth}
  		\centering
  		%
  		\captionof{figure}{\small$\vert P45,P3\vert\approx0.9896643320$, $\vert P45,P6\vert\approx1.0065187004$, $\vert P51,P50\vert\approx1.0146695273$}
  	\end{minipage}
  \end{center}
  {\par\centering  \par}

  \begin{center} 
  	\begin{minipage}{\linewidth}
      \centering
  	  %
  	  \captionof{figure}{\small$\vert P44,P40\vert\approx1.0912339657$, $\vert P47,P16\vert\approx1.0142532443$, $\vert P50,P49\vert\approx1.0908924633$}
  	\end{minipage}
  \end{center}
  {\par\centering  \par}

  \begin{center} 
  	\begin{minipage}{\linewidth}
  	  \centering
      %
      \captionof{figure}{52 vertices, asymmetric}
  	  \end{minipage}
  \end{center}
  {\par\centering$\vert P4,P54\vert\approx0.9029654473$, $\vert P38,P54\vert\approx1.0775714256$, $\vert P42,P41\vert\approx1.0445005488$\par}

  \begin{center} 
  	\begin{minipage}{\linewidth}
  		\centering
  		%
  		\captionof{figure}{\small$\vert P27,P52\vert\approx1.0818208359$, $\vert P27,P53\vert\approx1.0818208359$}
  	\end{minipage}
  \end{center}
  {\par\centering  \par}

  \begin{center} 
  	\begin{minipage}{\linewidth}
  		\centering		
  		%
  		\captionof{figure}{\small$\vert P19,P53\vert\approx0.9903987194$, $\vert P44,P54\vert\approx0.9903987194$}
  		
  		
  	\end{minipage}
  \end{center}
  {\par\centering  \par}

  \begin{center} 
  	\begin{minipage}{\linewidth}
  		\centering	
  		\begin{tikzpicture}[draw=black,font=\sffamily\tiny,
  		myscale/.style={scale=#1, node distance=#1*5em},myscale=1]
  		\begin{scope}[xshift=4cm, myscale=1.3]
  		
  		\foreach \i/\x/\y in {
  			1/0.00000000000000000000/2.32325008948610411963,
  			2/0.54027144270551863414/1.48175928846870896827,
  			3/0.99888813208473958838/2.37039348329965449125,
  			4/1.53915957479025822252/1.52890268228225933989,
  			5/1.08054288541103726828/0.64026848745131381690,
  			6/1.99777626416947917676/2.41753687711320486287,
  			7/0.99172828738558715944/2.45160508159490442281,
  			8/0.38470545982402132301/3.24628947606805695614,
  			9/1.37643374720960864899/3.37464446817685725932,
  			10/0.76941091964804264602/4.16932886265001023673,
  			11/1.98345657477117431888/2.57996007370370428191,
  			12/1.75385243994376760490/1.37962918682269397586,
  			13/2.05750281089284703384/0.42684565830087578542,
  			14/2.73081236542557670433/1.16620635767225611090,
  			15/3.03446273637465635531/0.21342282915043775393,
  			16/3.70777229090738646988/0.95278352852181802390,
  			17/4.01142266185646612087/0.00000000000000000000,
  			18/2.69943185183290390583/1.70502060621078221914,
  			19/3.66845000738402315932/1.95201010844176225945,
  			20/2.96566124921244433921/2.66893030528581887140,
  			21/1.40079310577042170927/3.39385704817360656094,
  			22/1.75668030396461083420/4.32838596809074260818,
  			23/2.38806249008698978642/3.55291415361433982056,
  			24/2.74394968828117891135/4.48744307353147675599,
  			25/3.37533187440355808562/3.71197125905507308019,
  			26/3.73121907259774721055/4.64650017897221001562,
  			27/7.74264173445421288733/2.32325008948610500781,
  			28/7.20237029174869380910/3.16474089050349993713,
  			29/6.74375360236947418713/2.27610669567255508028,
  			30/6.20348215966395422072/3.11759749668994956551,
  			31/6.66209884904317117815/4.00623169152089797507,
  			32/5.74486547028473637511/2.22896330185900337639,
  			33/6.75091344706862628300/2.19489509737730603689,
  			34/7.35793627463019195289/1.40021070290415128312,
  			35/6.36620798724460534856/1.27185571079534942562,
  			36/6.97323081480616924210/0.47717131632220122217,
  			37/5.75918515968303967867/2.06654010526850484553,
  			38/5.98878929451044506038/3.26687099214951581772,
  			39/5.68513892356137162665/4.21965452067133295344,
  			40/5.01182936902863485074/3.48029382129995390471,
  			41/4.70817899807955519975/4.43307734982177148453,
  			42/4.03486944354682819380/3.69371665045039243580,
  			43/5.04320988262131120194/2.94147957276142690830,
  			44/4.07419172707018883983/2.69449007053044686799,
  			45/4.77698048524176854812/1.97756987368639114422,
  			46/6.34184862868379362055/1.25264313079859990196,
  			47/5.98596143048960271926/0.31811421088146746294,
  			48/5.35457924436722354500/1.09358602535786952892,
  			49/4.99869204617303442006/0.15905710544073373147,
  			50/4.36730986005065613398/0.93452891991713638031,
  			51/2.39607688924114681228/3.49086318385686444898,
  			52/5.34656484521306829549/1.15563699511534356823,
  			53/3.07619533581263970845/2.75776094089673051712,
  			54/4.66644639864157539932/1.88873923807547794418}
  		\coordinate (p-\i) at (\x,\y);
  		
  		\foreach \i/\j in {
  			2/1, 
  			3/1, 3/2, 
  			4/3, 4/2, 
  			5/4, 5/2, 
  			6/3, 6/4, 
  			7/1, 
  			8/1, 8/7, 
  			9/8, 9/7, 
  			10/8, 10/9, 
  			11/9, 11/7, 11/20, 
  			12/5, 
  			13/12, 13/5, 
  			14/12, 14/13, 
  			15/14, 15/13, 
  			16/14, 16/15, 
  			17/16, 17/15, 17/49, 17/50, 
  			18/6, 18/12, 
  			19/18, 19/16, 19/54, 
  			20/6, 20/18, 
  			21/10, 
  			22/10, 22/21, 
  			23/22, 23/21, 
  			24/22, 24/23, 
  			25/24, 25/23, 
  			26/24, 26/25, 26/41, 26/42, 
  			28/27, 
  			29/27, 29/28, 
  			30/28, 30/29, 
  			31/28, 31/30, 
  			32/29, 32/30, 
  			33/27, 
  			34/27, 34/33, 
  			35/33, 35/34, 
  			36/34, 36/35, 
  			37/33, 37/35, 
  			38/31, 
  			39/31, 39/38, 
  			40/38, 40/39, 
  			41/39, 41/40, 
  			42/40, 42/41, 
  			43/32, 43/38, 
  			44/42, 44/43, 44/53, 44/54, 
  			45/32, 45/43, 45/52, 45/37, 
  			46/36, 
  			47/36, 47/46, 
  			48/46, 48/47, 
  			49/47, 49/48, 
  			50/48, 50/49, 
  			51/21, 51/11, 51/20, 
  			52/46, 52/37, 
  			53/25, 53/51, 53/19, 
  			54/50, 54/52}
  		\draw[gray,thin] (p-\i) -- (p-\j);
  		
  		\draw[red,thin] (p-11) -- (p-20);
  		\draw[red,thin] (p-45) -- (p-37);
  		\foreach \i in {1,...,54}
  		\fill[black] (p-\i) circle (1.1pt);
  		
  		\foreach \i/\a in {
  			1/153,
  			2/153,
  			3/93,
  			4/333,
  			5/198,
  			6/33,
  			7/337,
  			8/157,
  			9/337,
  			10/159,
  			11/337,
  			12/78,
  			13/318,
  			14/78,
  			15/318,
  			16/78,
  			17/219,
  			18/285,
  			19/343,
  			20/45,
  			21/219,
  			22/39,
  			23/339,
  			24/39,
  			25/279,
  			26/138,
  			27/333,
  			28/333,
  			29/333,
  			30/153,
  			31/93,
  			32/345,
  			33/97,
  			34/37,
  			35/277,
  			36/339,
  			37/157,
  			38/258,
  			39/18,
  			40/198,
  			41/78,
  			42/258,
  			43/105,
  			44/163,
  			45/225,
  			46/39,
  			47/219,
  			48/159,
  			49/219,
  			50/159,
  			51/95,
  			52/275,
  			53/285,
  			54/105}
  		\node[anchor=\a] (P\i) at (p-\i) {\i};
  		\end{scope}
  		\end{tikzpicture}
  		\captionof{figure}{\small$\vert P11,P20\vert\approx0.9862260008$, $\vert P45,P37\vert\approx0.9862260008$}
  		
  		
  	\end{minipage}
  \end{center}
  {\par\centering  \par}

  \begin{center} 
  	\begin{minipage}{\linewidth}
  	  \centering

  	  \captionof{figure}{\small$\vert P49,P54\vert\approx1.0587135610$, $\vert P50,P53\vert\approx1.0587135610$}
  	\end{minipage}
  \end{center}
  {\par\centering  \par}

  \begin{center} 
  	\begin{minipage}{\linewidth}
  	  \centering
  	  %
  	  \captionof{figure}{\small$\vert P13,P40\vert\approx1.0664925618$, $\vert P27,P53\vert\approx1.0664925618$}
  	\end{minipage}
  \end{center}
  {\par\centering  \par}

   \begin{center} 
  	\begin{minipage}{\linewidth}
  	  \centering
  	  %
  	  \captionof{figure}{\small$\vert P19,P52\vert\approx1.0897807877$, $\vert P51,P44\vert\approx1.0897807877$}
  	\end{minipage}
  \end{center}
  {\par\centering  \par}

  \begin{center} 
  	\begin{minipage}{\linewidth}
  	  \centering
  	  %
  	  \captionof{figure}{\small$\vert P25,P55\vert\approx1.0142467355$, $\vert P26,P54\vert\approx1.0142467355$, $\vert P49,P56\vert\approx1.0142467355$, $\vert P50,P53\vert\approx1.0142467355$}
  	\end{minipage}
  \end{center}
  {\par\centering  \par}

  \begin{center} 
  	\begin{minipage}{\linewidth}
  	  \centering
  	  %
  	  \captionof{figure}{\small$\vert P46,P57\vert\approx1.0126823415$, $\vert P55,P44\vert\approx1.0126823415$, $\vert P56,P48\vert\approx1.0126823415$}
  	\end{minipage}
  \end{center}
  {\par\centering  \par}

  \begin{center} 
  	\begin{minipage}{\linewidth}
  	  \centering
  	  %
  	  \captionof{figure}{\small$\vert P50,P57\vert\approx1.0510811050$, $\vert P51,P55\vert\approx1.0510811050$, $\vert P56,P49\vert\approx1.0510811050$}
  	\end{minipage}
  \end{center}
  {\par\centering  \par}

  \begin{center} 
  	\begin{minipage}{\linewidth}
  	  \centering
  	  %
      \captionof{figure}{\small$\vert P56,P58\vert\approx1.0319785178$, $\vert P57,P29\vert\approx1.0319785178$, $\vert P59,P57\vert\approx1.1635893883$, $\vert P59,P58\vert\approx1.1635893883$} 
  	\end{minipage}
  \end{center}
  {\par\centering  \par}

  \begin{center} 
  	\begin{minipage}{\linewidth}
  	  \centering
  	  %
  	  \captionof{figure}{59 vertices, asymmetric}
  	\end{minipage}
  \end{center}
  {\par\centering$\vert P50,P59\vert\approx1.0143213484$, $\vert P56,P51\vert\approx1.0184958217$, $\vert P56,P58\vert\approx1.0273592144$\par}

  \begin{center}
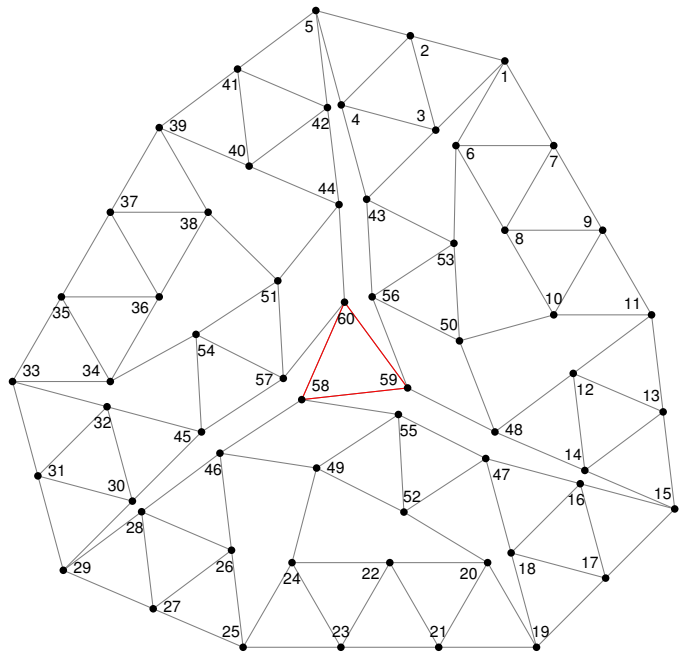
 
  	\begin{minipage}{\linewidth}
  	  \centering
  	  %
  	  \captionof{figure}{60 vertices, rotational symmetry of order 3}
  	\end{minipage}
  \end{center}
  {\par\centering$\vert P58,P60\vert=\vert P58,P59\vert=\vert P59,P60\vert\approx1.0889437642$\par}

  \begin{center} 
  	\begin{minipage}{\linewidth}
  	  \centering
  	  %
  	  \captionof{figure}{60 vertices, point symmetry}
  	\end{minipage}
  \end{center}
  {\par\centering$\vert P29,P59\vert=\vert P30,P58\vert\approx0.9013368132$\par}

  \begin{center} 
  	\begin{minipage}{\linewidth}
  	  \centering
  	  %
  	  \captionof{figure}{61 vertices, mirror symmetry}
  	\end{minipage}
  \end{center}
  {\par\centering$\vert P31,P60\vert\approx0.8891556455$, $\vert P31,P61\vert=\vert P60,P61\vert\approx1.0097449640$\par}

  \begin{center} 
  	\begin{minipage}{\linewidth}
  	  \centering
  	  %
  	  \captionof{figure}{61 vertices, point symmetry}
  	\end{minipage}
  \end{center}
  {\par\centering$\vert P27,P29\vert=\vert P56,P58\vert\approx1.0396024985$, $\vert P60,P61\vert\approx0.9805271527$\par}

  \begin{center} 
  	\begin{minipage}{\linewidth}
  		\centering
  		%
  		\captionof{figure}{62 vertices, point symmetry}
  	\end{minipage}
  \end{center}
  {\par\centering$\vert P59,P62\vert=\vert P61,P60\vert\approx1.3001527612$\par}

  \begin{center} 
  	\begin{minipage}{\linewidth}
  	  \centering
  	  %
  	  \captionof{figure}{62 vertices, point symmetry}
  	\end{minipage}
  \end{center}
  {\par\centering$\vert P25,P62\vert=\vert P56,P32\vert\approx1.0332092374$\par}

  \begin{center} 
  	\begin{minipage}{\linewidth}
  	  \centering
  	  %
  	  \captionof{figure}{62 vertices, mirror symmetry}
  	\end{minipage}
  \end{center}
  {\par\centering$\vert P61,P27\vert=\vert P62,P56\vert\approx1.0055442035$, $\vert P61,P62\vert\approx1.3434011015$\par}

  \begin{center} 
  	\begin{minipage}{\linewidth}
  	  \centering
  	  %
  	  \captionof{figure}{62 vertices, mirror symmetry} 	  
  	\end{minipage}
  \end{center}
  {\par\centering$\vert P28,P61\vert=\vert P61,P56\vert\approx1.0758636209$, $\vert P30,P60\vert=\vert P60,P58\vert\approx0.9760901087$\par}

%


  \newpage
  
  \section{How to find minimal 4-regular matchstick graphs}
  
  We are interested in all new solutions on 4-regular matchstick graphs. If you find a new graph, a better approximate solution for an already existing graph, or a proof, please submit it to the author's institutional E-Mail address or website contact.
  \\ \\
  Here we give a brief instruction for the interested reader how minimal 4-regular matchstick graphs could be found.
  
  \begin{itemize}
  	
  	\item Take a close look at all graphs found so far \cite{Catalogue}. Study the design from the outer circle to the center. Identify rigid and flexible subgraphs (triangles, rhombuses, etc.).
  	
  	\item First construct the frame of the graph using only equilateral triangles. Then built bigger rigid subgraphs of 1, 3, 5, or 7 triangles. Construct the graph from outside to inside. It is always helpful using a mirror or point symmetry for the frame. Asymmetric graphs are the most difficult types to handle.
  	
  	\item Use the free software \textsc{Matchstick Graphs Calculator} (MGC) \cite{MGC} available on the author's website. The MGC contains a brief description of the construction language and a manual for using the function buttons.
  	
  	\item Use a CAD software. All graphs from this article and many others graphs are included into the MGC and can be downloaded as DXF, SVG, and some other files.
  	
  	\item Experiment with already existing graphs. Use subgraphs for creating new graphs.
  	
    \item Do not use matchsticks. Even they give the name of such graphs, matchsticks are the worst possible method to construct such graphs. They are simply too inaccurate.
    
  	\item Better use flexi filing strip fasteners. Some of the author's graphs have been found with this kind of model. Pictures can be found in \cite{Thread2}.
  	
  	\item Be creative and develop new methods. Perhaps automated methods using artificial intelligence are the key to success.
  	
  \end{itemize}
  
  \newpage
  
  \section{References}
  
  \begingroup
  \renewcommand{\section}[2]{}
  
  \endgroup
  
\end{document}